\documentclass[12pt,a4paper,final]{amsart}
\usepackage{polski}
\usepackage{stmaryrd}
\usepackage{amssymb,amsmath}
\usepackage{amscd}
\usepackage{enumerate}
\usepackage[OT4]{fontenc}
\usepackage{graphicx}
\input xy
\xyoption{all}

\def\rzym{\begin{enumerate}[i)]}
\def\krzym{\end{enumerate}}
\def\rzymm{\begin{enumerate}[(i)]}
\def\krzymm{\end{enumerate}}
\def\ara{\begin{enumerate}[1)]}
\def\kara{\end{enumerate}}
\def\araa{\begin{enumerate}[(1)]}
\def\karaa{\end{enumerate}}
\def\lit{\begin{enumerate}[a)]}
\def\klit{\end{enumerate}}
\def\litt{\begin{enumerate}[(a)]}
\def\klitt{\end{enumerate}}
\def\cent{\begin{center}}
\def\kcent{\end{center}}

\def\symbolindex#1#2{}

\newtheorem {proposition}{Proposition}
\newtheorem{tw}[proposition]{Theorem}
\newtheorem{wn}[proposition]{Corollary}
\newtheorem{cw}[proposition]{Exercise}
\newtheorem{lem}[proposition]{Lemma}
\newtheorem{df}[proposition]{Definition}

\theoremstyle{definition}
\newtheorem*{ex}{Example}
\newtheorem*{exs}{Examples}
\theoremstyle{remark}
\newtheorem*{rem}{Remark}
\newtheorem*{rems}{Remarks}

\def\dfi{\begin{df}}
\def\kdfi{\end{df}} 
\def\twi{\begin{tw}}
\def\ktwi{\end{tw}}
\def\stw{\begin{proposition}}
\def\kstw{\end{proposition}}
\def\wni{\begin{wn}}
\def\kwni{\end{wn}}
\def\cwi{\begin{cw}}
\def\kcwi{\end{cw}}
\def\uwa{\begin{rem}}
\def\kuwa{\end{rem}}
\def\uwi{\begin{rems}}
\def\kuwi{\end{rems}}
\def\prz{\begin{ex}}
\def\kprz{\end{ex}}
\def\pry{\begin{exs}}
\def\kpry{\end{exs}}
\def\dowo{\begin{proof}}
\def\kdowo{\end{proof}}
\def\strz{\rightarrow}
\def\equa{\begin{equation}}
\def\equab{\begin{equation*}}
\def\kequa{\end{equation}}
\def\kequab{\end{equation*}}
\def\lemm{\begin{lem}}
\def\klemm{\end{lem}}

\newcounter{num}
\newenvironment{cid}[1] {\minCDarrowwidth#1pt $\begin{CD}}{\end{CD}$\ }
\def\cda#1{\begin{cid}{#1} }
\def\cd{\begin{cid}{10} }
\def\kcd{\end{cid}}

\makeindex

\xyoption{rotate}

\title{Zredukowane homologie Khovanova}

\author[W. Lubawski]{Wojciech Lubawski}
\address{Theoretical Computer Science Department\\
Jagiellonian University\\ Go³êbia 24\\ 00-300 Kraków, Poland}
\address{Institute of Mathematics\\ Polish Academy of Sciences\\
Œniadeckich 8 \\ 00-956 Warszawa, Poland}
\email{w.lubawski@gmail.com}

\date{}
\begin{document}

\begin{abstract}
W³aœciwie od momentu zdefiniownia przez Khovanova \cite{khov1} nowego niezmiennika bêd¹cego kategoryfikacj¹ wielomianu Jonesa, nazwanego póŸniej homologiami Khovanova, zauwa¿ono potrzebê (i to zarówno teoretyczn¹ jak i obliczeniow¹) pos³ugiwania siê \emph{zredukowanymi homologiami}. Niestety, okaza³o siê ¿e jedyne dostêpne definicje wprzypadku zwyk³ych homologii s¹ nienaturalne -- zale¿¹ od dokonywania wyborów sk³adowych odpowiednich diagramów i nie daj¹ siê uogólniæ na struktury kategoryjne wy¿szych poziomów (jak na przyk³ad na morfizmy miêdzy wêz³ami); zawód ten spotêgowa³ fakt skonstruowania takich zredukowanych homologii dla nieparzystych homologii Khovanova podanych przez Rasmussena, Ozsvatha i Szabo w \cite{ROS}. W poni¿szej pracy konstuujemy \emph{w³aœciwe} zredukowane homologie Khovanova, które mo¿na naturalnie zdefiniowaæ dla dowolnego diagramu wêz³a lub splotu.   
\end{abstract}



\maketitle

\section{Wprowadzenie}

W \cite{Put} zdefiniowano now¹ wersjê homologii Khovanova. Przypomnijmy skonstruowane tam struktury algebraiczne.

Rozwa¿my specjalny rodzaj algebry tensorowej generowanej przez modu³ z gradacj¹ $A = R\mathbf{1} \varoplus R\mathbf{x}$ nad $R=\mathbb Z[X,Y,Z] / X^2=Y^2=1$ z operacjami danymi poprzez:
\begin{itemize}
\item mno¿enie $m\colon A^{\varotimes 2} \strz A^{\varotimes 1}$
\begin{align*} m(\mathbf{1}\mathbf{1}) &= 1 & m(\mathbf{1}\mathbf{x}) &= \mathbf{x}\\
m(\mathbf{x}\mathbf{1})&= XZ\mathbf{x} & m(\mathbf{x}\mathbf{x}) &= 0  \end{align*} 

\item komno¿enie $\Delta\colon A^{\varotimes 1}\strz A^{\varotimes 2}$
\begin{align*} 
\Delta (\mathbf{1}) &= \mathbf{x}\mathbf{1} + YZ \mathbf{1}\mathbf{x}\\
\Delta(\mathbf{x}) &= (\mathbf{x}\mathbf{x})\end{align*} 

\item jednoϾ $\eta\colon R\strz A$
\begin{align*}
\eta(1)&=\mathbf{1}
\end{align*}
\item kojednoϾ $\epsilon\colon A\strz R$
\begin{align*}
\epsilon(\mathbf{1})&=0&\epsilon(\mathbf{x}) &=1 \end{align*}
\item permutacja $P\colon A^{\varotimes 2}\strz A^{\varotimes 2}$
\begin{align*}
P(\mathbf{1}\mathbf{1}) &= X\mathbf{1}\mathbf{1} & P(\mathbf{1}\mathbf{x}) &= Z^{-1}\mathbf{x}\mathbf{1} \\
P(\mathbf{x}\mathbf{x}) &= Y\mathbf{x}\mathbf{x} & P(\mathbf{x}\mathbf{1}) &= Z \mathbf{1}\mathbf{x}
\end{align*}
\end{itemize}

Oznaczmy przez $D$ dowolny diagram rozpatrywanego splotu. Diagram taki wyposa¿amy w ka¿dym skrzy¿owaniu w strza³ki (mo¿emy wybraæ jeden z dwóch nierównowa¿nych sposobów narysowania takich strza³ek, niemniej oba sposoby dadz¹ izomorficzne obiekty) zgodnie z nastêpuj¹c¹ zasad¹:

\includegraphics{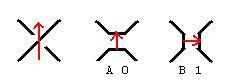}

Obrazek powy¿ej pokazuje równie¿ sposoby tzw. wyg³adzenia skrzy¿owania. Jeœli mamy na diagramie $D$ dok³adnie $n$ skrzy¿owañ, to indeksuj¹c dowolnie skrzy¿owania otrzymamy $n$ strza³ek przypisanych do $n$ skrzy¿owañ w diagramie $D$. Podobnie dla dowolnego $I\subseteq \{ 0,1 \}^n$ otrzymamy $n$ strza³ek na $n$ skrzy¿owaniach, przy czym $k$-te skrzy¿owanie wyg³adzamy w sposób wskazany przez $k$-t¹ wspó³rzêdn¹ $I$. Wyg³adzenie takie oznaczymy $D(I)$, sk³ada siê ono z dok³adnie $k(I)$ roz³¹cznych okrêgów (które bêdziemy dalej oznaczaæ $\overline{s}_1,\ldots ,\overline{s}_{k(I)}$). Dla u³atwienia opisu przyjmijmy, ¿e $s(a_i)$ oznacza okr¹g z którego $a_i$ wychodzi, $t(a_i)$ okr¹g do którego $a_i$ idzie.

Niech teraz $K_{k}$ oznacza dla ka¿dego $k$ graf pe³ny o $k$ wierzcho³kach, wyposa¿ony w dodatkowe informacje. Jeœli $k=k(I)$ to u¿ywamy tak¿e oznaczenia $K_I$.

Zauwa¿my na przyk³ad, ¿e $D(0)$ mo¿emy uto¿samiæ z podgrafem $K_{k(0)}$ wyposa¿onym dodatkowo w cykliczny porz¹dek pomiêdzy krawêdziami przy ka¿dym z wierzcho³ków. Podgraf ten bêdziemy oznaczaæ $\mathcal M(D(0))$.

Niech $\mathcal A'(I)$ bêdzie zbiorem sk³adaj¹cym siê z wierzcho³ków i krawêdzi $K(I)$. Rozpinaj¹c na tym zbiorze abelow¹ grupê woln¹, a nastêpnie na grupie jej pierœcieñ grupowy otrzymamy pierœcieñ $\mathcal A(I)$, który bêdziemy dalej nazywaæ woln¹ algebr¹ strza³ek $D(I)$. Niech $A$ bêdzie pierœcieniem Frobeniusa po raz pierwszy zdefiniowanym przez Khovanova w \cite{}, tj. $A= \mathbb Z \mathbf{1}\varoplus \mathbb Z \mathbf{x}$ gdzie $A$ jest wyposa¿one w dodatkow¹ gradacjê, przy której $\deg \mathbf{1} = -1$ oraz $\deg \mathbf{x} = +1$. Definiujemy reprezentacjê $\mathcal A$ do odpowiedniej potêgi tensorowej $A$ w sposób nastêpuj¹cy:
$$ ev_I\colon \mathcal A\strz \mathbb Z Mor (K^{\varotimes k(I)}, K^{\varotimes k(I)})$$
Jeœli $a_i\in\mathcal A$ jest strza³k¹ dla której $s(a_i) = \overline{s}_s$ oraz $t(a_i) = \overline{s}_t$ (dla $s\neq t$) to przyjmujemy $ev_I(a_i) = T^{s,t}$; jeœli $a_i$ taka, ¿e $s(a_i) = t(a_i) = \overline{s}_s$ to $ev_I(a_i) = T^s$ gdzie
$T^{s,t}\colon K^{\varotimes 2}\strz K^{\varotimes 2}$ jest dane poprzez
\begin{align*}
T^{s,t}(\mathbf{1}\mathbf{1}) &= \mathbf{x}\mathbf{1} + \mathbf{1}\mathbf{x} & T^{s,t}(\mathbf{1}\mathbf{x}) &= \mathbf{x}\mathbf{x} \\
T^{s,t}(\mathbf{x}\mathbf{1}) &= \mathbf{x}\mathbf{x} & T^{s,t}(\mathbf{x}\mathbf{x}) &= 0 \mathbf{1}\mathbf{x}
\end{align*}
ponadto $T^{s}\colon K\strz K$ jest dane jako
\begin{align*}
T^{s}(\mathbf{1}) &= 2\mathbf{x}& T^{s}(\mathbf{x}) &= 0
\end{align*} 
Rozszerzamy $ev_I$ do ca³ego $\mathcal A$ przyjmuj¹c $ev_I (a_{i_1}\ldots a_{i_m}) = ev_I(a_{i_1})\circ \ldots \circ ev_I(a_{i_m})$. Oznaczmy poprzez $\mathcal O_I$ obraz odwzorowania $ev_I$.

Zauwa¿my, ¿e mamy nastêpuj¹ce operacje. Jeœli $D(I)$ oraz $D(J)$ s¹ dwoma wyg³adzeniami diagramu $D$ ró¿ni¹cymi siê tylko sposobem wyg³adzenia skrzy¿owania o indeksie $i$, to definiujemy $\partial^I_J\colon \mathcal A(I)\strz \mathcal B(I) $ jako:
$$ \partial_J^I (a_{i_1} \ldots a_{i_m}) = \begin{cases} a_{i_1} \ldots a_{i_m} & ,s(a_i)\neq t(a_i) \text{ w }I\\ a_i a_{i_1} \ldots a_{i_m} & ,s(a_i) = t(a_i) \text{ w } I \end{cases} $$

Zauwa¿my, ¿e zachodzi nastêpuj¹ce

\twi
Nastêpuj¹cy diagram jest przemienny

\includegraphics{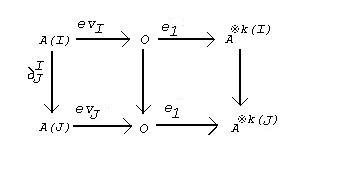}

gdzie $e_{1}$ jest wziêciem wartoœci danego odwzorowania na odpowiednio du¿ej potêdze tensorowej jedynki, a ró¿niczka z prawej strony jest ró¿niczk¹ Khovanova.  
\ktwi

Przyjmuj¹c jako $|I|$ sumê wszystkich wspó³czynników w $I$ oraz $$\partial^I = \sum_{I<J} \epsilon_J^I \partial^I_J,$$ gdzie $\epsilon^I_J$ jest uznakowieniem Khovanova otrzymamy kompleks $$\mathcal O_k\colon = \oplus_{|I|=k}\mathcal O_I.$$

\section{Dowód niezmienniczoœci -- kompleks strza³kowy}

W celu wykazania, ¿e kompleks $\mathcal O$ jest niezmienniczy wzglêdem ruchów Reidemaistera przedstawimy alternatywny opis tego kompleksu w jêzyku grafów. Podobny zabieg zastosowa³ Bloom w pracy \cite{bloom} dowodz¹c niezmienniczoœci homologii nieparzystych wzglêdem mutacji. Niech zatem $\mathcal R'(D(I))$ bêdzie zbiorem pografów $K_I$ zdefiniowanym w sposób nastêpuj¹cy. $X\in \mathcal R(D(I)) \Leftrightarrow$ gdy dla dowolnej sk³adowej spójnej $S$ zachodzi jeden z poni¿szych warunków:
\ara
\item $S$ jest wyró¿nionym wierzcho³kiem;
\item $\pi(S) = 0$ i $S$ ma co najwy¿ej jeden wyró¿niony wierzcho³ek;
\item $\pi(S) = \mathbb Z$ i $S$ nie posiada wyró¿nionych wierzcho³ków.
\kara
Definiujemy $\mathcal R(D(I))\colon = \mathbb Z\mathcal R'(D(I))$
Zauwa¿my, ¿e mamy odwzorowanie 
$$\psi\colon \mathcal R(D(I))\strz \mathcal A(I)$$
przyporz¹dkowuj¹ce strza³ce $a_i$ odwzorowanie $T^{s(a_i),t(a_i)}$ lub $T^{s(a_i)=t(a_i)}$, wierzcho³kowi $\overline{s}_m$ odwzorowanie $T^{m}$, a podgrafowi sk³adaj¹cemu siê z paru krawêdzi i wierzcho³ków odpowiednie z³o¿enie (zauwa¿my, ¿e kolejnoœæ z³o¿enia nie jest istotna).

W wielu miejscach poni¿ej bêdziemy rysowaæ fragmenty diagramów w którch -- bia³y wierzcho³ek oznacza niewyró¿niony wierzcho³ek diagramu, czarny wierzcho³ek oznacza wyró¿niony wierzcho³ek, przerywana strza³ka oznacza strza³kê $D(I)$ nie le¿¹c¹ w rozpatrywanym jego podzbiorze, a ci¹g³a strza³ka oznacza strza³kê w $D(I)$ która le¿y w rozpatrywanym jego podzbiorze. 

\lemm
J¹dro odwzorowania $\psi$ jest generowane przez nastêpuj¹ce relacje (przedstawione lokalnie):

\includegraphics{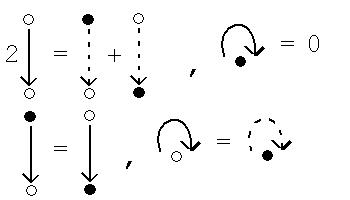}
\newline dodatkowo jeœli $l_1 , \ldots , l_{2m}$ jest parzystym cyklem w grafie, to $l_1 + l_3 + \ldots + l_{2m-1} = l_2 +\ldots +l_{2m}$, a jeœli $l_1 , \ldots , l_{2m}$ jest dowolnym cyklem w grafie, to $l_1 + l_2 + \ldots + l_{m} = l_1 +\ldots +l_{2m-1} +$ dowolny wyró¿niony wierzcho³ek na cyklu (tak¿e rozumiane jako relacje lokalne).
\klemm

\dfi
Zbiór $\Theta_I\colon = \mathcal R(D(I)) / \ker \psi$
nazwiemy grup¹ grafów dla diagramu $D(I)$.
\kdfi

Zauwa¿my ¿e mamy, podobnie jak wczeœniej, ró¿niczkowanie. Jeœli wyg³adzenia diagramu $D(I)$ oraz $D(J)$ ró¿ni¹ siê tylko wyg³adzeniem w skrzy¿owaniu $i$-tym (lub inaczej, wzd³u¿ strza³ki $a_i\in\mathcal M(D(0))$) to w zale¿noœci od sytuacji wyjœciowej mo¿emy mieæ dwie sytuacje:

\includegraphics{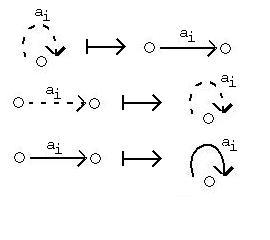}
\newline podobnie jeœli któryœ z przylegaj¹cych wierzcho³ków jest wyró¿niony to wyró¿niamy tak¹ sam¹ liczbê wierzcho³ków po drugiej stronie (zauwa¿my ¿e zgodnie z opisem j¹dra $\psi$ ró¿niczka jest poprawnie okreœlona).

PrzejdŸmy nastêpnie do wykazania, ¿e kompleks (a dok³adniej jego homologie) s¹ rzeczywiœcie niezmiennikiem wêz³a. Dowód przeprowadzimy w duchu Bar-Natana \cite{barnatan}, korzystaj¹c silnie z jego intuicji geometrycznej.

\subsection{Niezmienniczoœæ wzglêdem 1RM}

W przypadku pierwszego ruchu Reidemaistera mamy nastêpne odwzorowanie kompleksów:

\includegraphics{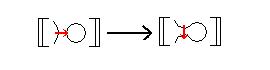}
\newline Warto zauwa¿yæ, ¿e nasz wyjœciowy kompleks jest sto¿kiem nad tym odwzorowaniem.

odwzorowanie zdefiniowane w sposób nastêpuj¹cy:

\includegraphics{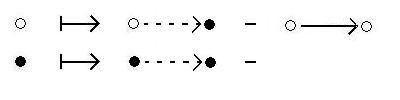}
\newline indukuje izomorfizm homologii szukanego wêz³a i wêz³a ,,bez pêtelki'', tj. po wykonaniu pierwszego ruchu Reidemaistera.

\subsection{Niezmienniczoœæ wzglêdem 2RM}

Mamy nastêpuj¹cy diagram:

\includegraphics{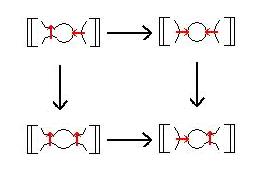}
\newline niestety, w naszym przypadku mamy dwie mo¿liwoœci uzupe³nienia tego diagramu do pe³nych okrêgów. Niezmienniczoœæ musimy pokazaæ dla ka¿dego z nich zwracaj¹c uwagê na przemiennoœæ tych dwóch przypadków z odpowiednimi ró¿niczkami.  

Obydwa przypadki przedstawiamy na nastêpuj¹cych dwóch diagramach:

\begin{center}
\includegraphics{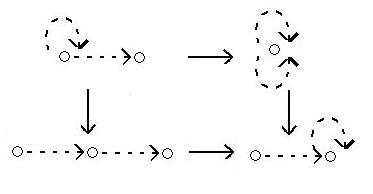}
\newline
\includegraphics{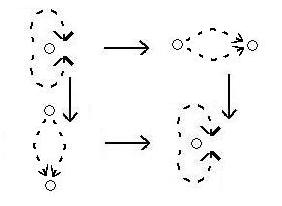}
\end{center}
a odpowiednie odwzorowania, prowadz¹ce do przek¹tnej diagramów (kompleks lewy dolny jest w parze pierwszy)
przedstawiamy jako:

\begin{center}
\includegraphics{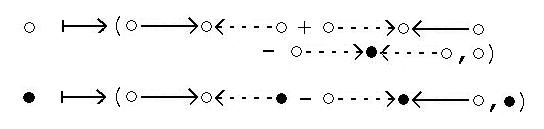}
\newline
\includegraphics{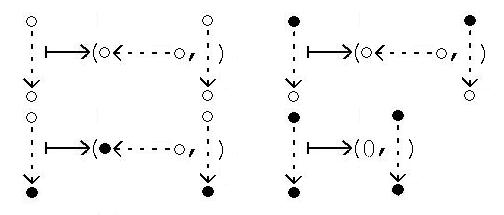}
\end{center}

Z tego, ¿e pomiêdzy dowolnymi dwoma okrêgami jest co najmniej jedna strza³ka oraz z postaci j¹dra $\psi$ wynika, ¿e odwzorowania te sk³adaj¹ siê w jedno odwzorowanie kompleksów, indukuj¹ce izomorfizmy w homologiach, co dowodzi niezmienniczoœci.

\subsection{Niezmienniczoœæ wzglêdem 3RM}

Zauwa¿my wpierw, ¿e odwzorowanie zdefiniowane dla 2RM jest w rzeczywistoœci retrakcj¹ deformacyjn¹, wiêc stosuj¹c twierdzenie o sto¿kach otrzymamy, ¿e kompleksy dla dwóch stron trzeciego ruchu Reidemaistera s¹ homotopijne do sto¿ków:

\includegraphics{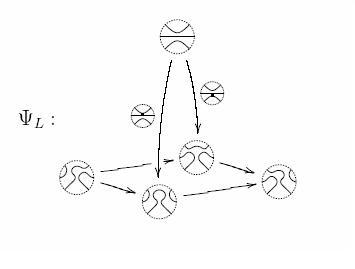}
\newline
\includegraphics{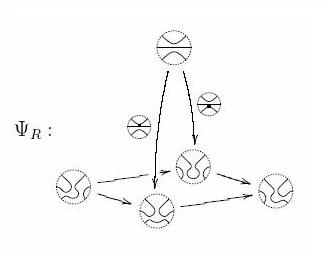}
\newline a te, jak ³atwo widaæ, s¹ jednakowe.

\section{Wersja chronologiczna}

Niech $R' = \mathbb Z [X,Y]/ X^2 = Y^2 = 1$. Zgodnie z powy¿szym podrozdzia³em wystarczy rozpatrywaæ kompleksy nad $R'$ zamiast nad $R$. £atwo zauwa¿yæ (czego ju¿ nie bêdziemy robiæ), ¿e wszystkie powy¿sze definicje i rozumowania da siê przeprowadziæ równie¿ w przypadku chronologicznym. Istotna mo¿e byæ równie¿ informacja, ¿e taki zredukowany chronologiczny kompleks Khovanova w przypadku nieparzystych homologii Khovanova redukuje siê do zredukowanych homologii Khovanova (zdefiniowanych oryginalnie w \cite{ROS}).

\section{Wêz³y wirtualne}

Wêz³y wirtualne to wêz³y wyposa¿one w jeszcze jeden rodzaj skrzy¿owania, tzw. skrzy¿owanie wirtualne. Dwa diagramy wêz³a wirtualnego s¹ równowa¿ne, jeœli od jednego do drugiego mo¿emy dojœæ poprzez ruch Reidemaistera oraz tzw. detour move: fragment wêz³a pomiêdzy dwoma punktami, który posiada wy³¹cznie skrzy¿owania wirtualne mo¿emy wymazaæ i narysowaæ w dowolny sposób ponownie - pamiêtaj¹c by wszystkie nowe skrzy¿owania oznaczyæ jako wirtualne.

Homologie dla wêz³ów wirtualnych definiujemy poprzez wziêcie dwóch kopii kompleksu zredukowanego (drug¹ kopiê bêdziemy oznaczaæ poprzez dodanie kreski pod jednym z wierzcho³ków) oraz dodefiniowanie nowego odwzorowania, tzw. twistu id¹cego od jednego okrêgu do jednego okrêgu jako dodanie kreski do danego elementu zgodnie z relacjami (kreskê mo¿emy zawsze dowolnie ustawiæ na diagramie): 

\includegraphics{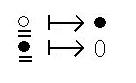}

Aby wykazaæ niezmienniczoœæ homologii dla wêz³ów wirtualnych ³atwo zauwa¿yæ, ¿e wystarczy sprawdziæ 2RM. Istotnie, do pierwszego ruchu wirtualnoœæ wêz³a nic nie wnosi, a trzeci wynika w taki sam sposób z drugiego jak w przypadku klasycznym. Drugi ruch w jedynym nieklasycznym przypadku wygl¹da w sposób nastêpuj¹cy:

\includegraphics{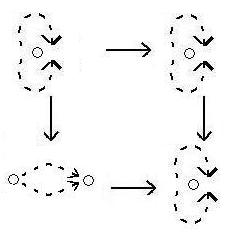}
\newline Licz¹c homologie dochodzimy, podobnie jak i w pozosta³ych przypadkach, posi³kuj¹c siê interpretacj¹ geometryczn¹ nastêpuj¹ce odwzorowanie:

\includegraphics{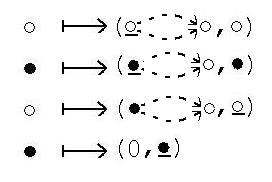}
\newline £atwo tak¿e zauwa¿yæ, ¿e jest ono zgodne z odwzorowaniami zdefiniowanymi poprzednio czyli sk³ada siê z nimi do odwzorowania kompleksów, a to dowodzi niezmienniczoœci wzglêdem 2RM.


\begin{thebibliography}{Wi2'}
\bibitem{barnatan} D.~Bar-Natan, \emph{Khovanov's homology for tangles and cobordisms}, Geom. Topol. 9 (2005), 1443-1499, arXiv:math/0410495v2
\bibitem{khov1} M.~Khovanov, \emph{A categoryfication of the Jones polynomial}, Duke Math. J. 101 (2000), 359-426, arXiv:math/9908171v2
\bibitem{khov2} M.~Khovanov, \emph{An invariant of tangle cobordisms}, Trans. Amer. Math.
Soc. 358 (2006), 315-327,

\bibitem{ROS} P.~Ozvath, J.~Rasmussen, Z.~Szabo, \emph{Odd Khovanov homology}, arXiv:0710.4300v1
\bibitem{Put} K.~Putyra, \emph{Cobordisms with chronologies and a generalisation of the Khovanov complex}, arXiv:1004.0889v1
\bibitem{shum1} A.~Shumakovitch, \emph{Patterns in odd Khovanov homology}, arXiv:1101.5607v2
\bibitem{carSai1} J.S.~Carter, M.~Saito, \emph{Knotted surfaces and their diagrams}, Math. Surv. Mon. 55, AMS, 1998

\bibitem{cMW1} D.~Clark, S. Morrison and K. Walker, \emph{Fixing the functoriality of Khovanov
homology}, arXiv:math.GT/0701339.

\bibitem{caprau1} C.~Caprau, \emph{An sl(2) tangle homology and seamed cobordisms}, arXiv:0707.3051.

\bibitem{man1} V.O.~Manturov, \emph{Khovanov homology for virtual knots with
arbitrary coefficients}, arXiv:math/0601152v3

\bibitem{krmro} P.B.~Kronheimer, T.S.~Mrowka, \emph{Khovanov homology is an unknot-detector}, arXiv:1005.4346 

\bibitem{bloom} J.~Bloom, \emph{Odd Khovanov homology is mutation invariant}, arXiv:0903.3746

\bibitem{tavares} D.~Tavares, \emph{Odd Khovanov homology for virtual knots}, dostêpne online: https://dspace.ist.utl.pt/bitstream/2295/722646/1/Dissertacao.pdf

\bibitem{turtur1} V.~Turaev, P.~Turner, \emph{Link homology and unoriented topological quantum field theory},
Algebr. Geom. Topol. 6 (2006), 1069.1093


\end{thebibliography}
\end{document}